\documentclass[12pt]{amsart}
\numberwithin{equation}{section}
\usepackage{amssymb}
\usepackage{latexsym}
\def\R{\mathbb R}
\def\C{\mathbb C}

\def\N{\mathbb N}
\def\FF{{\mathcal F}}

\def\im{\operatorname{Im}}

\def\meas{\operatorname{meas}}
\newtheorem{la}{Lemma}
\numberwithin{la}{section}
\newtheorem*{tha}{Theorem A}
\newtheorem*{thb}{Theorem B}

\newtheorem{thm}{Theorem}

\theoremstyle{remark}

\begin{document}
\title{Fixed points of composite entire and quasiregular maps}
\author{Walter Bergweiler}
\thanks{
Part of this research was 
done while the author was visiting Purdue University.
He thanks Alexandre Eremenko and the Purdue Department
of Mathematics for the hospitality.
Support
by the Alexander von Humboldt Foundation, the National Science
Foundation (grant of A.~Eremenko),
and the G.I.F.,
the German--Israeli Foundation for Scientific Research and
Development, Grant G -809-234.6/2003, is gratefully acknowledged.}
\address{Mathematisches Seminar,
Christian--Albrechts--Universit\"at zu Kiel,
Ludewig--Meyn--Str.\ 4,
D--24098 Kiel,
Germany}
\curraddr{ Department of Mathematics, Purdue University,
150 N.\ University Street, West Lafayette, IN 47907, USA}
\email{bergweiler@math.uni-kiel.de}
\begin{abstract}
We give a new proof of the result that if
$f$ and $g$ are entire transcendental functions, then 
$f\circ g$ has infinitely many fixed points.
The method yields a number of generalizations of this result. 
In  particular, it extends to quasiregular maps 
in $\R^d$.
\end{abstract}
\maketitle

\section{Introduction and main results}
The following result was conjectured by Gross 
(see~\cite[p.~542]{Ehr} and~\cite[Problem~5]{Gro72})
and first proved in~\cite{Ber90}.

\begin{tha}
Let $f$ and $g$ be entire transcendental functions.
Then $f\circ g$ has infinitely many fixed points.
\end{tha}

The following generalization of Theorem~A was proved in~\cite{Ber91}.
\begin{thb}
Let $f$ and $g$ be entire transcendental functions.
Then $f\circ g$ has infinitely many repelling fixed points.
\end{thb}
Here a fixed point $\xi$ of a holomorphic function $h$ is
called repelling if $|h'(\xi )|>1$. The repelling fixed points play
an important role in iteration theory.

The purpose of this paper is twofold. Firstly, we give a new proof of Theorem~A.
Secondly, we obtain some generalizations of Theorem~A (and~B).

The main difference between the method employed here  and the previous proofs 
of Theorem~A and~B is that 
the Wiman-Valiron method which was crucial 
in~\cite{Ber90,Ber91} is not used here.
Instead we use some ideas from normal families. 
This method is also applicable for quasiregular maps; see~\cite{Rick}
for the definition and basic properties of quasiregular maps.
\begin{thm} \label{quasireg}
Let $d\geq 2$ and let $f,g:\R^d\to \R^d$ be quasiregular maps with an essential 
singularity at $\infty$.
Then $f\circ g$ has infinitely many fixed points.
\end{thm}
For functions in the plane we also obtain some extensions of the previously known results.
\begin{thm} \label{largemult}
Let $f$ and $g$ be entire transcendental functions.
Then  there exists a sequence $(\xi_n)$ such that 
$(f\circ g)(\xi_n)=\xi_n$ and
$(f\circ g)'(\xi_n)\to\infty $.
\end{thm}
\begin{thm} \label{nonreal}
Let $f$ and $g$ be entire transcendental functions.
Then $f\circ g$ has infinitely many nonreal fixed points.
\end{thm}
Theorem~\ref{nonreal} answers a question of Clunie~\cite{Clu96} 
who had shown that at least one of the two functions 
$f\circ g$ and $g\circ f$ has infinitely many nonreal fixed points.
The special case $f=g$ had been dealt with earlier in~\cite{BCL},
answering a question of Baker~\cite{Bak66}.
Theorem~\ref{nonreal} implies that for any straight line
there are infinitely many fixed points not lying on this line.

Similar ideas to the ones employed in this paper were used -- in the context
of iteration rather than composition -- in~\cite{BB,BerQ,Ess98,Ess00}
for holomorphic maps, and in~\cite{Sie,Sie2} for quasiregular maps.

Although each of the Theorems~\ref{quasireg}--\ref{nonreal} contains Theorem~A as a
special case, we will first give a proof of Theorem~A in~\S\ref{proofA}, as this
explains the underlying idea best. In \S\S \ref{proof1}--\ref{proof3} we will then
prove Theorems~\ref{quasireg}--\ref{nonreal}. These sections will make occasional reference
to~\S\ref{proofA},  but are independent of each other.

\section{Proof of Theorem A}\label{proofA}
\subsection{Preliminary Lemmas}\label{ahlforssec}
We shall need a result from the Ahlfors theory of covering surfaces; 
see~\cite{Ahl35}, \cite[Chapter 5]{Hay64} or \cite[Chapter XIII]{Nev53}
for an account of this theory.
To state the result of the Ahlfors theory that we need,
let $D\subset \C$ be a domain and let $f:D\to\C$ be holomorphic.
Given a Jordan domain $V\subset\C$, we say that $f$ has 
an {\em island} over $V$ if $f^{-1}(V)$ has a component whose closure
is contained in $D$. Note that if $U$ is such a component, 
then 
$f|_U:U\to V$ is a proper map.
\begin{la} \label{ahlfors} 
Let $D\subset \C$ be a domain and let $D_1,D_2\subset\C$
be Jordan domains
with disjoint closures.
Let $\FF$ be a family of functions holomorphic in $D$ 
which is not normal.
Then there exists a function $f\in\FF$ which has an island over
$D_1$ or $D_2$.
\end{la}
For example,
Lemma \ref{ahlfors} follows from
Theorem~5.5 (applied with a domain $D_3$ 
containing~$\infty$) and 
Theorem~6.6 in~\cite{Hay64}.

For a different proof 
of Lemma~\ref{ahlfors}
see~\cite[\S 5.1]{Ber98}.
The proof given there is particularly simple in the case where
 the $D_j$ are small 
disks. It turns out that this special case suffices for our purposes.

The following lemma is a simple consequence of the maximum  principle.
\begin{la} \label{quasi} 
Let $D\subset \C$ be a domain and let 
$(f_n)$ be sequence of functions holomorphic in $D$ which is not normal.
If $(f_n)$ converges locally uniformly in $D\backslash E$ for some finite set $E$,
then $f_n\to\infty$ in $D\backslash E$.
\end{la}
This lemma  will be useful when dealing with  quasinormal families.
By definition, a family $\FF$ of functions holomorphic in a domain
$D$ is called {\em quasinormal} 
(cf.~\cite{Chu,Mon,Schi})
if 
for each sequence $(f_n)$ in $\FF$ there exists a 
subsequence $(f_{n_k})$ and
a finite set $E\subset D$ such that
$\left(f_{n_k}\right)$ converges locally uniformly in
$D\backslash E$.
If the cardinality of the exceptional set $E$ can be bounded
independently of the sequence $(f_n)$, and
if $q$ is the smallest 
such bound,
then we
say that $\FF$ is quasinormal of {\em order}~$q$.

We denote the maximum modulus of an entire function $f$ by $M(r,f)$.
\begin{la} \label{MAr} 
Let $f$ be an entire transcendental function and $A>1$. Then
$$\lim_{r\to\infty}\frac{M(Ar,f)}{M(r,f)}=\infty.$$
\end{la}
This result follows easily 
from the convexity of
$\log M(r,f)$ in $\log r$ and the transcendency of $f$.
We omit the details. For an alternative proof of Lemma~\ref{MAr} 
see the proof of Lemma~\ref{MArqr} in \S\ref{proof1lemmas} below.

\subsection{Proof of Theorem A}
We first choose a sequence $(c_n)$ tending to $\infty$ such that $|f(c_n)|\leq 1$.
We may assume that $|c_n|\geq |g(0)|$ for all $n$ 
and define $r_n$ by $M(r_n,g)=|c_n|$.

We define
$$
f_n(z):=\frac{f(c_n z)}{r_n}\ \ \ \text{and}\ \ \ 
g_n(z):=\frac{g(r_n z)}{c_n}.$$ 
It is easy to see that no subsequence of $(f_n)$
is normal at $0$.
Since $f_n(1)\to 0$ it follows from Lemma~\ref{quasi} that 
$(f_n)$ is not normal in $\C\backslash\{0\}$.
Passing to a subsequence if necessary we may thus assume that no 
subsequence of $(f_n)$ is normal at $a_1:=0$ and some $a_2\in\C\backslash\{0\}$.

It follows from Lemma~\ref{MAr} that if $n\to\infty$, then 
$M(r,g_n)\to 0$ if $r<1$ and 
$M(r,g_n)\to \infty$ if $r>1$. 
Lemma~\ref{quasi} implies that the sequence $(g_n)$ is not quasinormal.
Passing to a subsequence if necessary we may thus assume
that there exist $b_1,b_2,b_3\in\C\backslash\{0\}$
where no subsequence of $(g_n)$ is normal.
We choose $0<\varepsilon<\frac12$
such that the 
closed disks of radius $\varepsilon$ around the $b_j$
are pairwise disjoint and do not contain~$0$.
In the following we denote by $B(a,r)$ the open disk of radius
$r$ around a point $a$; that is,
$B(a,r):=\{z\in \C:|z-a|<r\}$. 

It follows from Lemma~\ref{ahlfors} that if $n$ is sufficiently
large and $j\in\{1,2\}$, then $f_n$ has an island 
in $B(a_j,\varepsilon)$
over at least two of the domains $B(b_k,\varepsilon)$.
This implies 
that there exists  $k\in\{1,2,3\}$ such that 
$f_n$ has an island  $U_1\subset B(a_1,\varepsilon)$
and another  island  $U_2\subset B(a_2,\varepsilon)$
over the same disk $B(b_k,\varepsilon)$.

Moreover, 
it follows from Lemma~\ref{ahlfors} that if $n$ is sufficiently
large, then there exists
$j\in\{1,2\}$ such that $g_n$ has an island 
$V\subset B(b_k,\varepsilon)$ over
$B(a_j,\varepsilon)$.
Then $V\cap  g_n^{-1}(U_j)$ contains a component $W$ of 
$(f_n\circ g_n)^{-1}(B(b_k,\varepsilon))$
satisfying $\overline{W}\subset B(b_k,\varepsilon)$.

For $z\in\partial W$ we have 
$$|((f_n\circ g_n)(z)-b_k)-\left((f_n\circ g_n)(z)-z\right)|=
|z-b_k|<\varepsilon=|(f_n\circ g_n)(z)-b_k|.$$
Rouch\'e's theorem implies that the number of fixed points of
$f_n\circ g_n$ in $W$ coincides with the number of zeros of 
$f_n\circ g_n-b_k$ in $W$.
As $f_n\circ g_n$ is a proper map from $W$ onto $B(b_k,\varepsilon)$ the function
$f_n\circ g_n-b_k$ has at least one zero in $W$
and thus 
 $f_n\circ g_n$ has a fixed point
$\xi\in W\subset B(b_k,\varepsilon)$. Then $\xi_n:=r_n \xi$ is 
a fixed point of $f\circ g$. Since $\xi_n\in B(r_nb_k,r_n\varepsilon)$
it follows that $\xi_n\to\infty$ as $n\to\infty$
so that $f\circ g$ has infinitely
many fixed points.\qed

\section{Proof of Theorem~\ref{quasireg}}\label{proof1}
\subsection{Preliminary Lemmas} \label{proof1lemmas}
As a general reference for quasiregular maps 
we re\-commend~\cite{Rick}.
We first state some lemmas analogous to those stated in~\S\ref{ahlforssec}, and begin
with the analogue of Lemma~\ref{ahlfors}.

We note that Lemma~\ref{ahlfors} is a generalization of Montel's theorem,
which in turn is the result that corresponds to Picard's theorem in the context
of normal families.
The analogue of Picard's theorem for quasiregular maps was given by 
Rickman~\cite{Ric80} who proved
that there exists $q=q(d,K)\in\N$ with the property that every $K$-quasiregular
map $f:\R^d\to\R^d$ which omits $q$ points is constant.
We shall call this number $q$ the {\em Rickman constant}.
The corresponding normality result was proved by Miniowitz~\cite{Min}, using 
an extension of the Zalcman lemma~\cite{Zal75} to quasiregular maps.
We refer to~\cite{Min} also for further information about normal families of
quasiregular maps. Besides normal families 
we will also consider quasinormal families of quasiregular maps, which are
defined in exactly the same way as for holomorphic functions.

Miniowitz's extension of the Zalcman lemma 
has been used by Siebert~\cite{Sie,Sie2} to deduce the following Lemma~\ref{covqr}
from Rickman's theorem.
Here we call, as in~\S\ref{ahlforssec}, a domain $U$ an 
{\em island} of the
quasiregular map $f:D\to\R^d$ over the
simply connected domain $V\subset \R^d$,
if $U$ is a component of $f^{-1}(V)$ and if $\overline{U}\subset D$.
And as in dimension~$2$ we denote by  $B(a,r)$ the open ball of radius
$r$ around a point $a\in\R^d$; that is, $B(a,r):=\{x\in \R^d:|x-a|<r\}$. 
Here $|x|$ is the (Euclidean) norm of a point $x\in\R^d$.
With this notation Siebert's result (see~\cite[Satz 2.2.2]{Sie} or
\cite[Corollary~3.2.2]{Sie2}) can be stated as follows.
\begin{la} \label{covqr} 
Let $d\geq 2$, $K\geq 1$ and let $q=q(d,K)$ be the Rickman constant.
Let $a_1,\dots,a_q\in\R^d$ be distinct.
Then there exists $\varepsilon>0$ with the following property:
if  $D\subset \R^d$ is a domain 
and $\FF$ is a non-normal family of functions $K$-quasiregular in $D$, 
then there exists a function $f\in\FF$ which has an island over
$B(a_j,\varepsilon)$ for some $j\in\{1,\dots,q\}$.
\end{la}
The following lemmas is literally the same as Lemma~\ref{quasi} in \S\ref{ahlforssec}, 
and it is again a simple consequence of the maximum principle.
\begin{la} \label{quasiqr} 
Let $D\subset \R^d$ be a domain and let 
$(f_n)$ be non-normal sequence of functions 
which are $K$-quasiregular in $D$.
If $(f_n)$ converges locally uniformly in $D\backslash E$ for some finite set $E$,
then $f_n\to\infty$ in $D\backslash E$.
\end{la}
The next lemma is identical to Lemma~\ref{MAr} in \S\ref{ahlforssec}.
Again we denote the by $M(r,f)$ the maximum modulus; that is,
$M(r,f):=\max_{|x|=r}|f(x)|$.
\begin{la} \label{MArqr} 
Let $f:\R^d\to\R^d$ be quasiregular with an essential singularity 
at~$\infty$ and let $A>1$. Then
$$\lim_{r\to\infty}\frac{M(Ar,f)}{M(r,f)}=\infty.$$
\end{la}
\begin{proof}
Suppose that the conclusion does not hold. Then there exist $C>1$ and a sequence
$(r_n)$ tending to $\infty$ such that $M(Ar_n,f)\leq C M(r_n,f)$.
The sequence $(f_n)$ defined by
$$f_n(x):=\frac{f(r_n x)}{M(r_n,f)}$$
is then bounded and thus normal in $B(0,A)$. Passing to a subsequence we may
assume that $f_n\to h$ for some quasiregular map $h:B(0,A)\to \R^d$. 
We have $h(0)=0$ while $M(1,h)=1$. Thus $h$ is not constant. 

Now there exists $a\in\R^d$ such that $f$ has infinitely many $a$-points.
Without loss of generality we may assume that $a=0$ since otherwise
we can consider $f(x+a)-a$ instead of $f(x)$.
A contradiction will now be obtained from 
Hurwitz's theorem (cf.~\cite[Lemma~2]{Min}).

More precisely, choose $0<t<1$ such that $h(x)\neq 0$ for 
$|x|=t$.  
For sufficiently large $n$ we then have
$\mu(h,B(0,t),0)=\mu(f_n,B(0,t),0)
=\mu(f,B(0,r_n t),0)$. Here
$\mu(h,B(0,t),0)$ denotes the topological degree.
Thus
$$\mu(h,B(0,t),0)=\sum_{x\in h^{-1}(0)\cap B(0,r)}i(x,h)$$
where $i(x,h)$ is the topological index.
But $\mu(f,B(0,r_n t),0)\to\infty$ as $n\to\infty$ since
$f$ has infinitely many zeros. This is 
a contradiction.
\end{proof}

The next lemma is a simple consequence of Lemma~\ref{MArqr}.
\begin{la} \label{fastgrowth} 
Let $f:\R^d\to\R^d$ be quasiregular with an essential singularity 
at~$\infty$. Then
$$\lim_{r\to\infty}\frac{\log M(r,f)}{\log r}=\infty .$$
\end{la}

The following lemma
(see~\cite[Lemma~1.3.14]{Sie} or~\cite[Lemma~2.1.5]{Sie2}) 
replaces the argument where Rouch\'e's theorem was used in the proof
of Theorem~A. 
\begin{la} \label{rouchetype} 
Let $U\subset \R^d$ be a domain and let $a\in\R^d$ and $r>0$ be such
that  $\overline{U}\subset B(a,r)$.
Suppose that $h:U\to B(a,r)$ is proper and quasiregular.
Then $h$ has a fixed point in $U$.
\end{la}
Next we shall need the following result.
\begin{la} \label{minmax} 
Let $K>1$, let $D\subset \R^d$ a domain and let
$C$ be a compact subset of~$D$.
Then there exist $\alpha,\beta>0$ with the following property:
if $f$ is $K$-quasiregular in  $D$ 
and satisfies $|f(x)|\geq 1$ for all $x\in D$, 
then 
$\log|f(y)|\leq \alpha+\beta \log|f(x)|$ for all $x,y\in C$.
\end{la}
\begin{proof}
It follows from~\cite[Corollary 3.9, p.~91]{Rick} that there exist $A,B>0$ such 
that if $B(a,2\delta)\in D$,
then $\log|f(x)|\leq A+B\log|f(a)|$ for all
$x\in B(a,\delta)$.
We may assume that $C$ is connected.
Since $C$ is compact there exist $a_1,\dots,a_N\in C$  and $\delta>0$
with
$$C\subset \bigcup_{j=1}^N B(a_j,\frac12 \delta)\ \ \ \text{and} \ \ \  
\bigcup_{j=1}^N B(a_j,2\delta)\subset D.$$
The conclusion follows with $\beta:=B^{N+1}$ and some $\alpha$.
\end{proof}
Finally we need the following observation apparently
made first in~\cite[Lemma~3]{Gro68}.
\begin{la} \label{symm} 
Let $A,B$ be sets and let  $f:A\to B$ and $g:B\to A$ be functions.
Then the set of fixed points of $f\circ g$ and 
the set of fixed points of $g\circ f$ have the same cardinality.
\end{la}
To prove this lemma we only have to observe that $g$ is a bijection
from the set of fixed points of $f\circ g$ to 
the set of fixed points of $g\circ f$.

\subsection{Proof of Theorem~\ref{quasireg}}
Let $(c_n)$ be a sequence in $\R^d$ which
tends to $\infty$ 
and define
$F_n(x):=f(|c_n|x)/|c_n|$ and $G_n(x):=g(|c_n|x)/|c_n|$.
Lemma~\ref{fastgrowth} yields that $M(r,F_n)\to\infty$ and $M(r,G_n)\to\infty$ 
as $n\to\infty$ if $r>0$, while $F_n(0)\to 0$ and $G_n(0)\to 0$.
Thus
no subsequence of $(F_n)$ or $(G_n)$ is normal at~$0$. 

We distinguish between two cases.

{\em Case 1}. For every choice of  $(c_n)$ the
 sequences $(F_n)$ and $(G_n)$ are both quasinormal. 

We may choose the sequence $(c_n)$ such that $|g(c_n)|\leq 1$ for all $n$.
Applying Lemma~\ref{quasiqr} and, 
passing to subsequences if necessary, we may assume that 
$F_n\to\infty$ in $\R^d\backslash E_f$ and that
$G_n\to\infty$ in $\R^d\backslash E_g$ for two finite sets 
$E_f$ and $E_g$ containing $0$. Moreover, we may assume that
$E_g$ contains at
least one point $b\in\R^d$ with $|b|=1$,
with no subsequence of $(G_n)$ converging in a neighborhood of~$b$.

We choose $\varepsilon>0$ such that $2\varepsilon<|a-b|$ for
all $a\in E_g\backslash\{b\}$ and $2\varepsilon<|a|$ for
all $a\in E_f\backslash\{0\}$. For sufficiently large $n$ 
we then have $|G_n(x)|> 2$ for $x\in\partial B(b,\varepsilon)$ and
$|F_n(x)|> 2$ for $x\in\partial B(0,\varepsilon)$,
while 
$|G_n(y)|<1$ for some  $y\in B(b,\varepsilon)$ and
$|F_n(z)|<1$ for some  $z\in B(0,\varepsilon)$.
Since $B(0,\varepsilon)\cup B(b,\varepsilon)\subset B(0,2)$ this implies that  
$G_n$ has an island $V\subset B(b,\varepsilon)$ over 
$B(0,\varepsilon)$ while
$F_n$ has an island $U\subset B(0,\varepsilon)$ over 
$B(b,\varepsilon)$.
As in the proof of Theorem~A we find that
$V\cap G_n^{-1}(U)$ contains a component $W$ of 
$(F_n\circ G_n)^{-1}(B(b,\varepsilon))$
satisfying $\overline{W}\subset B(b,\varepsilon)$.
Lemma~\ref{rouchetype} now implies that $F_n\circ G_n$ has 
a fixed point in $\xi\in B(b,\varepsilon)$
and thus 
$f\circ g$ has a fixed point $\xi_n \in B(|c_n|b,|c_n|\varepsilon)$. 
It follows that $f\circ g$ has infinitely many fixed points.

{\em Case 2}. The sequence $(c_n)$ can be chosen such that
one of the sequences $(F_n)$ and $(G_n)$ is not quasinormal. 

Because of Lemma~\ref{symm} we may assume that the sequence $(F_n)$
is not quasinormal. 
Passing to a subsequence if necessary, we may in fact assume that
no subsequence of $(F_n)$ is quasinormal. 

As in the proof of Theorem~A we may assume that $|c_n|\geq |g(0)|$ and 
define $r_n$ by $M(r_n,g)=|c_n|$.
As there we also define
$f_n(x):=f(|c_n| x)/r_n=F_n(x)|c_n|/r_n$ and $g_n(x):=g(r_n x)/|c_n|.$ 
Again we find that no subsequence of $(f_n)$ is normal at $0$. 

We now show that $(f_n)$ is not quasinormal.
To do this we assume that  $(f_n)$ is quasinormal. Passing to a subsequence
we then may assume that 
$f_{n}\to\infty$ in $\R^d\backslash E$ for some finite set $E$. 
Let $C\subset \R^d\backslash E$ be a compact set containing $\partial B(0,r)$ for 
some $r>0$. Then there exists a domain $D\supset C$ such that 
$|f_n(x)|\geq 1$ for $x\in D$ if $n$ is large.
Lemma~\ref{minmax} yields that
\begin{equation}\label{minmax2}
\log M(r,f_n)\leq \alpha +\beta \log |f_n(x)|
\end{equation}
for $x\in C$ and large $n$.
On the other hand, Lemma~\ref{fastgrowth} implies that
$$\log M(r,f_n)=\log M(|c_n|r,f)-\log r_n\geq 4\beta \log (|c_n|r) -\log r_n$$ 
if $n$ is large.
Lemma~\ref{fastgrowth} also yields
$\beta\log|c_n|=\beta\log M(r_n,g)\geq \log r_n$ for large $n$.
Thus $\log M(r,f_n)\geq 3\beta \log |c_n| +4\beta\log r\geq 2\beta \log |c_n| $
for large $n$. We deduce from this and (\ref{minmax2}) that 
$$\log |f_n(x)|\geq \frac{\log M(r,f_n)-\alpha}{\beta}
\geq  2\log |c_n| -\frac{\alpha}{\beta}\geq \log|c_n|$$
for $x\in C$ and large $n$. It follows that 
$$\log |F_n(x)|=\log |f_n(x)|-\log|c_n|+\log r_n\geq \log r_n$$
for $x\in C$ and large $n$. Hence $F_{n}\to\infty$ in $\R^d\backslash E$,
contradicting the assumption that no subsequence of $(F_n)$ is quasinormal. 
Thus $(f_n)$ is not quasinormal.
Passing to a subsequence if necessary we may assume that no subsequence
of $(f_n)$ is quasinormal.

From Lemma~\ref{MArqr} we deduce that
$M(r,g_n)\to 0$ if $r<1$ and 
$M(r,g_n)\to \infty$ if $r>1$.
Lemma~\ref{quasiqr} now
implies that no subsequence of $(g_n)$ is quasinormal.

Let $K$ be such that $f$ and $g$ are 
$K$-quasiregular and define $p:=2q-1\geq q$ where
$q=q(d,K)$ is the Rickman constant.
Passing to subsequence if necessary 
we may assume
that there exist $a_1,\dots,a_p\in\R^d\backslash\{0\}$
where no subsequence of $(f_n)$ is normal and
that there exist $b_1,\dots,b_p\in\R^d\backslash\{0\}$
where no subsequence of $(g_n)$ is normal.

It follows from Lemma~\ref{covqr} that there exists 
$\varepsilon>0$ such that if 
$n$ is sufficiently large and $j\in\{1,\dots,p\}$, then $f_n$ has an island
in $B(a_j,\varepsilon)$ over 
at least $p-q+1$ of the $p$ balls $B(b_k,\varepsilon)$, 
and $g_n$ has an island in $B(b_j,\varepsilon)$ over 
at least $p-q+1$ of the $p$ balls $B(a_k,\varepsilon)$.

This implies that there exists $k\in\{1,\dots,p\}$ such that 
$f_n$ has an island in $B(a_j,\varepsilon)$ over $B(b_k,\varepsilon)$ for at 
least $p-q+1=q$ values of $j$. Lemma~\ref{covqr} implies that for at least one such value $j$ 
the function $g_n$ has an island in $B(b_k,\varepsilon)$ over $B(a_j,\varepsilon)$.

Thus we obtain $j,k\in\{1,\dots,p\}$ such that 
$f_n$ has an island
$U\subset B(a_j,\varepsilon)$ over $B(b_k,\varepsilon)$ and  
$g_n$ has an island
$V\subset B(b_k,\varepsilon)$ over $B(a_j,\varepsilon)$.
As 
before
we find that
$V\cap  g_n^{-1}(U)$ contains a component $W$ of 
$(f_n\circ g_n)^{-1}(B(b_k,\varepsilon))$
satisfying $\overline{W}\subset B(b_k,\varepsilon)$.
Lemma~\ref{rouchetype} now implies that $f_n\circ g_n$ has a fixed point
in $B(b_k,\varepsilon)$, and thus 
$f\circ g$ has a fixed point
in $B(|c_n| b_k,|c_n| \varepsilon)$.
Thus $f\circ g$ has infinitely many fixed points.\qed

\section{Proof of Theorem~\ref{largemult}}\label{proof2}
\subsection{Preliminary Lemmas}
We shall require an additional result from the Ahlfors theory.
An island $U$ of a function $f$ over a domain $V$ is called
{\em simple} if $f:U\to V$ is univalent. The following result  can also
be found in the references given in~\S\ref{ahlforssec}.
\begin{la} \label{ahlfors2}
Let $D\subset \C$ be a domain and let $D_1,D_2$ and
$D_3$ be Jordan domains
with pairwise disjoint closures.
Let $\FF$ be a non-normal family of functions holomorphic in $D$.
Then there exists a function $f\in\FF$ which has a simple island over
$D_1$, $D_2$ or $D_3$.
\end{la}
It follows from Lemma~\ref{ahlfors2} that a non-constant entire function $f$
has a simple island over one of three Jordan domains $D_1,D_2,D_3$
with pairwise disjoint closures. We shall need the simple observation made
in~\cite{BerQ} that we need need only two domains $D_1,D_2$ if $f$ is a 
polynomial or, more generally, a proper holomorphic map
whose range contains $D_1$ and $D_2$.

Although the formulation of the following Lemma~\ref{polylike}
was slightly different in~\cite[Lemma~2.2]{BerQ}, we omit
the simple proof based on the
Riemann-Hurwitz formula, but note that  it is analogous to the proof
of Lemma~\ref{polylike2}
which we will give below. 
\begin{la} \label{polylike}
Let $f:U\to V$ be a proper holomorphic map
and let $D_1$ and $D_2$ be Jordan domains
with disjoint closures contained in $V$.
Then there exist two domains $U_1,U_2\subset U$ which are simple
islands over $D_1$ or $D_2$.
\end{la}
Here $U_1$ and $U_2$ need not be islands over the same
domain. We allow the possibility that $U_1$ is an island over
$D_1$ and $U_2$ is an island over $D_2$, or vice versa. For example, this 
will always be the case if $f$ is univalent.
For proper maps of higher degree, however, we have 
the following lemma.
\begin{la} \label{polylike2}
Let $f:U\to V$ be a proper holomorphic map
of degree at least~$2$
and let $D_1$, $D_2$ and $D_3$ be Jordan domains
with pairwise disjoint closures contained in~$V$.
Then there exist  $k\in\{1,2,3\}$ such that $f$ has two simple
islands over $D_k$. 
\end{la}
\begin{proof}
Let $U_1,\dots,U_m$ be the components of
$f^{-1}\left(\bigcup_{k=1}^3 D_k\right)$. 
Thus the $U_j$ are the islands over the domains $D_k$.
Now $f|_{U_j}$ is a proper map of some degree $\mu_j$ and
$\sum_{j=1}^m \mu_j =3d$, where $d$ is the degree of $f$.
By the Riemann-Hurwitz formula the number of critical points
contained in $U_j$ is $\mu_j-1$, and $f$ has $d-1$ critical points in $U$.
Thus
$$3d-m=\sum_{j=1}^m (\mu_j -1)\leq d-1$$
so that $m\geq 2d+1$. Since $f$ has $d-1$ critical points in $U$
we conclude that the number $n$ of domains $U_j$ which do not contain a critical point
satisfies 
$$n\geq m-(d-1)\geq (2d+1)-(d-1)=d+2\geq 4.$$ 
Thus among the $U_j$ there at least $4$ simple
islands, and hence two of them must be over the same domain $D_k$.
\end{proof}

We shall also use the following 
well-known result; see, e.g.,~\cite[Lemma~2.3]{BerQ} for the simple proof.
\begin{la} \label{fix}
Let $0<\delta<\frac{\varepsilon}{2}$ and let 
$U\subset B(a,\delta)$ be a simply-connected domain.
Let $f:U\to B(a,\varepsilon)$  be holomorphic and bijective.
Then $f$ has a fixed point $\xi$ in $U$ which satisfies 
$|f'(\xi)|\geq\varepsilon/4\delta$.
\end{la}

We shall also need the following lemma concerning entire functions of small growth.
\begin{la} \label{smallgrowth}
Let $g$ be an entire function of the form 
$$g(z)=\prod_{k=1}^\infty \left( 1-\frac{z}{z_k}\right)$$
where $0<|z_1|\leq |z_2|\leq \dots$ and 
$\lim_{k\to\infty}|z_{k+1}/{z_k}|=\infty$.
Denote the zeros of $g'$ by $z_k'$, ordered such that
$0\leq |z_1'|\leq |z_2'|\leq \dots$. Then
$\lim_{k\to\infty}|z_{k+1}'/z_k'|=\infty$.
\end{la}
\begin{proof}
For sufficiently large $n$ there exists $r$ satisfying
$|z_n|\leq r$ and $8r\leq|z_{n+1}|$.
We show first that for such $r$ 
\begin{equation}\label{min}
\min_{|z|=4r}|g(z)|> \max_{|z|=r}|g(z)|,
\end{equation}
provided $n$ is large enough.
Let $|u|=r$ and $|v|=4r$. We will show that 
$|g(v)|>|g(u)|$. To this end  we write
\begin{eqnarray*}
\log\frac{|g(v)|}{|g(u)|}
&=&
\sum_{k=1}^\infty 
\log\frac{|v-z_k|}{|u-z_k|}\\
&=& \sum_{k=1}^n \log\frac{|v-z_k|}{|u-z_k|}
+ \sum_{k=n+1}^\infty \log\frac{|v-z_k|}{|u-z_k|}\\ &=&
S_1+S_2.
\end{eqnarray*}
For $k\leq n$ we have $|z_k|\leq r<|v|$ so that 
$$ \frac{|v-z_k|}{|u-z_k|} \geq \frac{|v|-|z_k|}{|u|+|z_k|} 
\geq
\frac{4r -r}{r +r}
= \frac{3}{2}.  $$
Thus $S_1\geq n\log\frac{3}{2}$.

For $k\geq n+1$ we have
$|z_k|\geq 8r>|v|$ 
so that 
$$
\log\frac{|v-z_k|}{|u-z_k|}
\geq
\log \frac{|z_k|-|v|}{|u|+|z_k|}
=
\log \left( 1 -\frac{4r}{|z_k|}\right)
-\log \left( 1 +\frac{r}{|z_k|}\right)
\geq 
-\frac{8r}{|z_k|}
-\frac{r}{|z_k|}.$$
Here we have used the inequalities $\log(1+x)\leq x$ and 
$\log(1-x)\geq -2x$ valid for $0\leq x\leq \frac12$.
For large $n$ we also have 
$|z_{k+1}/{z_k}|\geq 2$ if $k\geq n+1$. 
We find that
$$
\log\frac{|v-z_k|}{|u-z_k|}
\geq
-\frac{9r}{|z_k|}
\geq
-\frac{9r}{|z_{n+1}|}2^{n+1-k}
\geq
-\frac{9}{8}2^{n+1-k}
$$
for $k\geq n+1$. 
It follows that 
$$S_2\geq 
-\frac{9}{8}
\sum_{k=n+1}^\infty 2^{n+1-k} =
-\frac{9}{4}.$$
Together with the estimate for $S_1$ this implies that 
(\ref{min}) holds for large $n$.

Let now $U$ be the component of $g^{-1}(B(0,M(r,g)))$ which contains
$B(0,r)$. It follows from (\ref{min}) that $U\subset B(0,4r)$.
By our choice of $r$ the number of zeros of 
$g$ in $U$ is $n$.
The Riemann-Hurwitz formula yields that $g'$ has $n-1$ zeros in $U$.
Thus $g'$ has at most $n-1$ zeros in $B(0,r)$ and 
at least $n-1$ zeros in $B(0,4r)$. Thus $r\leq |z_n'|$ and
$|z'_{n-1}|\leq 4r$.
Since this holds for any $r$ satisfying 
$|z_n|\leq r$ and $8r\leq|z_{n+1}|$ we conclude,
choosing
$r=|z_n|$ or $r=\frac18 |z_{n+1}|$,
that 
$|z'_{n-1}|\leq 4|z_n|$
and $\frac18|z_{n+1}| \leq |z_n'|$ for large $n$.
Thus 
$\frac18|z_{n+1}| \leq |z_n'|\leq 4|z_{n+1}|$ for large $n$.
The conclusion follows.
\end{proof}

The following result is a 
variant of Lemma~\ref{symm}.
\begin{la} \label{gross2}
Let $f$ and $g$ be entire transcendental functions.
Suppose that there exists a sequence $(\xi_n)$ such that 
$(f\circ g)(\xi_n)=\xi_n$ and
$(f\circ g)'(\xi_n)\to\infty $.
Then $\eta_n:=g(\xi_n)$ satisfies 
$(g\circ f)(\eta_n)=\eta_n$ and
$(g\circ f)'(\eta_n)=(f\circ g)'(\xi_n)\to\infty $.
\end{la}
The proof is straight forward and thus omitted.

\subsection{Proof of Theorem~\ref{largemult}}
We proceed as in the proof of Theorem~A
and define the sequences $(c_n)$, $(r_n)$,
$(f_n)$ and $(g_n)$ as there.
Again we find that $(g_n)$ is not quasinormal in $B(0,2)$.
In the proof of Theorem~A we 
noted that by passing to a subsequence 
we can achieve that no subsequence of $(g_n)$ is normal at 
any of three points $b_1,b_2,b_3\in\C\backslash\{0\}$. The same argument
yields this for any number of points~$b_k$. Moreover, we can achieve
that these points are in $B(0,2)\backslash\{0\}$. 

We will have to distinguish several cases now. 

{\em Case 1}. It is possible to choose the sequence $(c_n)$ such 
that $(f_n)$ is not quasinormal of order $2$. 

Passing to a subsequence if necessary we may  assume
that there exist $a_1,a_2,a_3\in\C$
where no subsequence of $(f_n)$ is normal and 
that there exist $b_1,\dots,b_7\in\C\backslash\{0\}$
where no subsequence of $(g_n)$ is normal.

We choose $0<\varepsilon<\frac12$
such that the closures of the disks $B(a_j,\varepsilon)$ 
are pairwise disjoint. Moreover, we require 
that the closures of the disks $B(b_j,\varepsilon)$ are 
pairwise disjoint and do not contain $0$.
We also choose $0<\delta<\frac{\varepsilon}{2}$.

It follows from Lemma~\ref{ahlfors2} that if $n$ is sufficiently
large and $j\in\{1,2,3\}$, then $f_n$ has a simple island 
in $B(a_j,\varepsilon)$
over at least five of the
seven disks $B(b_k,\varepsilon)$.
Overall we obtain at least $15$ domains in the union of the 
three disks $B(a_j,\varepsilon)$
which are simple islands over one of the seven disks $B(b_k,\varepsilon)$.
This implies 
that there exists  $k\in\{1,\dots,7\}$ such that 
$f_n$ has three simple islands  $U_j\subset B(a_j,\varepsilon)$, $j\in\{1,2,3\}$, 
over the same disk $B(b_k,\varepsilon)$.
It also follows from Lemma~\ref{ahlfors} that if $n$ is sufficiently
large, then there exists
$j\in\{1,2,3\}$ such that $g_n$ has a simple island 
$V\subset B(b_k,\delta)$ over
$B(a_j,\varepsilon)$.
Then $W:=V\cap  g_n^{-1}(U_j)$ is a simple island of 
$f_n\circ g_n$ over $B(b_k,\varepsilon)$, 
and $W\subset B(b_k,\delta)$.
Lemma~\ref{fix} implies that $f_n\circ  g_n$ has a fixed point $\xi\in W$ with
$|(f_n\circ g_n)'(\xi)|\geq\varepsilon/4\delta$. Then $\xi_n :=r_n\xi$ is a fixed 
point of $f\circ g$ with 
$|(f\circ g)'(\xi_n)|\geq\varepsilon/4\delta$.
Since $\delta$ can be chosen arbitrarily small the conclusion follows.

{\em Case 2}.
For every choice of the sequence $(c_n)$ 
the sequence $(f_n)$ is quasinormal of order $2$. 

Then $f_n\to\infty$ in $\C\backslash\{0,1\}$ by  Lemma~\ref{quasi}. Let $K>2$.
For sufficiently large $n$ we then have $|f_n(z)|>2$ if
$\frac14\leq |z|\leq K$ and  $|z-1|\geq \frac14$.

We now distinguish two subcases.

{\em Case 2.1}. There are infinitely many $n$ such that
$B(1,\frac14)$ contains at least two zeros of $f_n$. 

Passing to a subsequence if necessary we may  assume
that this  holds for all $n$.
We note that if $n$ is sufficiently large, then $B(0,\frac14)$ also contains
at least two zeros of $f_n$. 
With $a_1:=0$ and $a_2:=1$ thus both disks $B(a_j,\frac14)$ contain at least two
zeros of $f_n$.

We may assume that no subsequence of $(g_n)$ is 
normal at five points $b_1,\dots,b_5\in B(0,2)\backslash\{0\}$. 
Again we choose $0<\varepsilon<\frac12$
such that the closures of the disks $B(b_k,\varepsilon)$ are 
pairwise disjoint and do not contain $0$, and we 
choose $0<\delta<\frac{\varepsilon}{2}$.
We can now deduce from Lemma~\ref{ahlfors}
that if $n$ is sufficiently large and $k\in\{1,\dots,5\}$,
then $g_n$ has an island over one of the disks $B(a_j,\frac14)$
in $B(b_k,\delta)$.
It follows that for at least three values of $k$ 
the function $g_n$ has an island over the 
same disk $B(a_j,\frac14)$
in $B(b_k,\delta)$. We may assume that the $b_k$ are numbered such that
this holds for $k\in\{1,2,3\}$.

We claim that there exist $k\in\{1,2,3\}$ such that $f_n$ has two simple 
islands $U_1,U_2$ in $B(a_j,\frac14)$ over  $B(b_k,\varepsilon)$. 
To this end we assume first that there are two components $X_1,X_2$ of $f_n^{-1}(B(0,2))$
contained in $B(a_j,\frac14)$ in which $f_n$ is univalent. Then we can simply take
$U_{1,2}:=f_n^{-1}(B(b_k,\varepsilon))\cap X_{1,2}$, for arbitrary $k\in\{1,2,3\}$.
Suppose now that such components $X_1,X_2$ do not exist.
Since $|f_n(z)|>2$ if
$|z-a_j|= \frac14$ and since $f_n$ has at least two zeros in $B(a_j,\frac14)$,
there now exists a component $X$ of $f_n^{-1}(B(0,2))$ contained in $B(a_j,\frac14)$ 
such that 
$f:X\to B(0,2)$ 
is a proper map of degree at least~$2$.
Now Lemma~\ref{polylike2} yields our claim
that there exists $k\in\{1,2,3\}$ such that 
$f_n$ has two simple 
islands $U_1,U_2\subset X\subset B(a_j,\frac14)$ over  $B(b_k,\varepsilon)$. 

Recall that in turn $g_n$ has an island $V\subset B(b_k,\delta)$ over 
$B(a_j,\frac14)$.  Lemma~\ref{polylike2}, applied to the proper map
$g_n:V\to B(a_j,\frac14)$,
now implies that $V$ contains a domain $W$ which is
a simple island of $g_n$ over $U_1$ or $U_2$. Then $W$ is 
a simple island of 
$f_n\circ g_n$ over $B(b_k,\varepsilon)$, 
and $W\subset V\subset B(b_k,\delta)$. As before
Lemma~\ref{fix} implies that $f_n\circ  g_n$ has a fixed point $\xi\in W$ with
$|(f_n\circ g_n)'(\xi)|\geq\varepsilon/4\delta$. 
Again $\xi_n :=r_n\xi$ is a fixed point of $f\circ g$ with
$|(f\circ g)'(\xi_n)|\geq\varepsilon/4\delta$, and the conclusion follows since
 $\delta$ can be chosen arbitrarily small. 

{\em Case 2.2}. For all sufficiently large $n$ 
the disk $B(1,\frac14)$ contains at most one  zero of $f_n$.

Since  $|f_n(z)|>2$ for
$\frac14\leq |z|\leq K$ and  $|z-1|\geq \frac14$, provided $n$ is sufficiently large,
we conclude that the annulus 
$R:=\{z\in\C: \frac14\leq |z|\leq K\}$ contains at most 
one zero of $f_n$. 
On the other hand, since $f_n(1)\to 0$ we conclude that $B(1,\frac14)$ 
contains a component $Y$ of $f_n^{-1}(B(0,2))$ and thus in
particular a zero of $f_n$ for large $n$. 
Thus $f_n$ has exactly one zero in 
$R$
if $n$ is sufficiently large.
In fact, $f_n$  takes every value 
in $B(0,2)$ exactly once in $Y$, 
and since there are no other 
components of $f_n^{-1}(B(0,2))$ intersecting 
$R$,
we see that
$f_n$  takes every value
in $B(0,2)$ exactly once in
$R$.
We find that if $z_1\in Y$ and $z_2\in \C$ with $|z_2|\geq |z_1|$ and
$f_n(z_2)=f_n(z_1)\in B(0,2)$, then $|z_2|\geq K$ while $|z_1|\leq
\frac54$ so that $|z_{2}|\geq \frac45 K|z_1|$.

In terms of $f$ we see
has infinitely many $c$-points for every 
$c\in\C\backslash\{0\}$, and the sequence $(w_n)$ of $c$-points
of $f$, arranged such that $|w_{n+1}|\geq |w_n|$, satisfies
$|w_{n+1}|\geq \frac45 K|w_n|$.
Since $K$ can be chosen arbitrarily large we deduce that
$\lim_{n\to\infty}|w_{n+1}/ w_n|=\infty$. 
It follows that the number $n(r,1/(f-c))$ 
of $c$-points of $f$ in $\overline{B(0,r)}$ satisfies
$n(r,1/(f-c))=o(\log r)$ as $r\to\infty$,
for any $c\in\C\backslash\{0\}$.
Standard estimates from value 
distribution theory~\cite{Hay64,Jan85,Nev53} 
now imply that $\log M(r,f)=o((\log r)^2)$
as $r\to\infty$. In particular, $f$ has order $0$.

Lemma~\ref{gross2} says that the conclusion is symmetric with 
respect to $f$ and $g$.
We may thus assume that 
interchanging the roles of $f$ and $g$ leads again to Case~2.2.
We find that
$g$ is also of order $0$ and that
the sequence $(z_n)$ of zeros of $g$ satisfies 
$\liminf_{n\to\infty}|z_{n+1}/ z_n|=\infty$. 

Assuming
without loss of generality that $g(0)\neq 0$ we 
deduce from Lemma~\ref{smallgrowth}
that the sequence $(z_n')$ of zeros of $g'$ satisfies
\begin{equation} \label{zn'}
\lim_{n\to\infty}\frac{|z_{n+1}'|}{|z_n'|}=\infty.
\end{equation}

Proceeding as before 
we may assume, passing to a subsequence if necessary, 
that no subsequence of $(g_n)$ is normal
at four  points $b_1,\dots,b_{4}\in B(0,2)\backslash \{0\}$.
As before 
we choose $0<\varepsilon<\frac12$
such that the closures of the disks $B(b_j,\varepsilon)$ are 
pairwise disjoint and do not contain $0$,
and we choose $0<\delta<\frac{\varepsilon}{2}$.

With $a_1:=0$ and $a_2:=1$ we can deduce as before  from Lemma~\ref{ahlfors}
that if $n$ is sufficiently large and $k\in\{1,\dots,4\}$,
then $g_n$ has an island 
$V_k\subset B(b_k,\delta)$
over one of the disks $B(a_j,\frac14)$.
It follows from (\ref{zn'}) that
at most one of these four islands 
contains a zero of $g_n'$, provided $n$ is large enough.
Thus $V_k$ is a simple island for 
at least three values of $k$.
Hence there exist two values of $k$ such that $V_k$ is a simple island 
over the same disk $B(a_j,\frac14)$. 
We may assume that the $b_k$ are numbered such that
this holds for $k\in\{1,2\}$.
By Lemma~\ref{polylike} the function $f_n$ has a simple island $U\subset B(a_j,\frac14)$
over $B(b_1,\varepsilon)$ or $B(b_2,\varepsilon)$.
Without loss of generality we can assume that $U$ is a simple island 
over $B(b_1,\varepsilon)$.
Then
$W:=V_1\cap  g_n^{-1}(U)$ is a simple island of 
$f_n\circ g_n$ over $B(b_1,\varepsilon)$, 
and $W\subset B(b_1,\delta)$.
As before 
we deduce that $f_n\circ g_n$ has a fixed point $\xi\in W$ with
$|(f_n\circ g_n)'(\xi)|\geq \varepsilon/4\delta$.
Again $\xi_n :=r_n\xi$ is a fixed point of $f\circ g$ with
$|(f\circ g)'(\xi_n)|\geq\varepsilon/4\delta$
so that the conclusion follows also in this case.\qed

\section{Proof of Theorem~\ref{nonreal}}\label{proof3}
\subsection{Preliminary Lemmas}
We shall need the following quantitative version of Lemma~\ref{ahlfors}.
\begin{la} \label{ahlfors3} 
Let $D_1,D_2\subset\C$ be Jordan domains with disjoint closures
and let $f:B(a,r)\to\C$ be a holomorphic function which has no
island over $D_1$ or $D_2$.
Then
$$
\frac{|f'(a)|}{2\mu(\log\mu+A)}
\leq
\frac{1}{r}
$$
where $\mu =\max\{1,|f(a)|\}$ and $A$ is a constant depending
only on the domains
$D_1$ and $D_2$.
\end{la}
While Lemma~\ref{ahlfors} says that the family of holomorphic
functions having no islands over any of two given domains 
is normal, Lemma~\ref{ahlfors3}
is based on the fact
that this family is in  fact a normal {\em invariant}
family.
Lemma~\ref{ahlfors3} follows from Lemma~\ref{ahlfors} together with results of
Hayman~\cite{Hay55} on normal invariant families.
It is a direct consequence
of Theorems~6.8, 6.6, and 5.5 of his book~\cite{Hay64}.

We also need the following version of a classical growth lemma due to Borel.
\begin{la} \label{borel} 
Let $r_0>0$ and let $T:[r_0,\infty)\to [e,\infty)$ be 
increasing and continuous.
Define 
$$F:=\left\{ r\geq r_0: T\left(r\left(1+\frac{1}{(\log T(r))^2}\right)\right)
>e T(r)\right\}.$$
Then 
$$\int_F \frac{dt}{t}\leq \frac{\pi^2}{6}.$$
\end{la}
\begin{proof}
We may assume that $F\neq \emptyset$ and define 
$$r_1:=\inf F, \ \ \ r_1':=r_1\left(1+\frac{1}{(\log T(r_1))^2}\right)$$
and then inductively
$$r_k:=\inf \left(F\cap [r_{k-1}',\infty)\right),
\ \ \ r_k':=r_k\left(1+\frac{1}{(\log T(r_k))^2}\right).$$
If $F\cap [r_{k}',\infty)=\emptyset$ for some $k$ so that the process terminates, then we put $N:=k$.
Otherwise we put $N:=\infty$.

We have
$T(r_k)\geq T(r_{k-1}')\geq  e T(r_{k-1})$. Hence
$T(r_k)\geq e^{k-1}T(r_1)\geq e^k$ so that 
$(\log T(r_k))^2\geq k^2$.
If $N=\infty$ we thus have $r_k\to\infty$ as $k\to\infty$.
In any case we find that
$$F\subset\bigcup_{k=1}^N [r_k,r_k']$$
so that
$$\int_F \frac{dt}{t}\leq 
\sum_{k=1}^N \log\frac{r_k'}{r_k}
=
\sum_{k=1}^N \log\left(1+\frac{1}{(\log T(r_k))^2}\right).$$
Since $\log(1+x)<x$ for $x>0$ this yields
$$\int_F \frac{dt}{t}\leq 
\sum_{k=1}^N \frac{1}{(\log T(r_k))^2}
\leq
\sum_{k=1}^\infty \frac{1}{k^2}=\frac{\pi^2}{6}.$$
\end{proof}

\subsection{Proof of Theorem~\ref{nonreal}}
We proceed as in the proof of Theorem~A and define the sequences 
$(c_n)$, $(r_n)$,
$(f_n)$ and $(g_n)$ and the points  $a_1,a_2$ and
$b_1,b_2,b_3$ as there. 
If we could choose all three points $b_1,b_2,b_3$ nonreal, then the argument
given there would imply that 
$f\circ g$ has infinitely many nonreal fixed points.
We may thus assume that $(g_n)$ is quasinormal of order $2$ in 
$\C\backslash\R$ and thus $g_n\to 0$ locally uniformly in 
$\C\backslash\R$ by Lemma~\ref{quasi}.

If there exist three 
annuli $\Omega_j:=\{z\in\C: S_j<|z|<T_j\}$  with disjoint closures
such that $g_n$ has an island over $B(a_1,\varepsilon)$ 
or $B(a_2,\varepsilon)$ in $\Omega_j\backslash \R$, for all $j\in\{1,2,3\}$
and infinitely many $n$, then the argument used in the proof of Theorem
A shows again that $f_n\circ g_n$ has a fixed point in one of these islands,
and thus $f\circ g$ has infinitely many nonreal fixed points.

We may thus assume that such annuli $\Omega_j$ do not exist. 
Passing to a subsequence if necessary we thus find 
an annulus $\Omega:=\{z\in\C: S-1<|z|<T+1\}$  with 
$S>1$ and $T>e^2 S$ such that 
no $g_n$ has an island over $B(a_1,\varepsilon)$ 
or $B(a_2,\varepsilon)$ in $\Omega\backslash \R$.
Since $g_n\to 0$ in $\C\backslash\R$ we may assume that
$|g_n(z)|\leq 1$ for $S\leq |z|\leq T$ and $|\im z|\geq 1$.

To save indices, we now write $h:=g_n$.
For $S\leq |z|\leq T$ and $|\im z|\leq 1$
we conclude from Lemma~\ref{ahlfors3} that
$$
\frac{|h'(z)|}{2\mu(\log\mu+A)}
\leq
\frac{1}{|\im z|}
$$
where $\mu =\max\{1,|h(z)|\}$ and $A$ is a constant.
We may assume that $A>1$.
If $S\leq r\leq T$ and $\log |h(re^{it})|\geq A$,  then 
$|h(re^{it})|\geq e^A>1$ so that
$|\im(re^{it})|=r|\sin t|\leq 1$ and thus
\begin{equation}\label{hlogh}
\frac{|h'(re^{it})|}{|h(re^{it})|\log|h(re^{it})|}
\leq
\frac{2|h'(re^{it})|}{|h(re^{it})|(\log|h(re^{it})|+A)}
\leq
\frac{4}{r|\sin t|}.
\end{equation}
We denote by $T(r,h)$ the Nevanlinna characteristic of $h$
and recall the inequality
\begin{equation}\label{MT}
\log M(r,h)\leq \frac{R+r}{R-r}T(R,h)
\end{equation}
valid for $0<r<R$.
Since $M(r,g_n)\to \infty$ as $n\to\infty$ if $r>1$
we conclude that 
$T(r,h)=T(r,g_n)\to\infty$ for $r>1$.
In particular we may assume that $T(r,h)\geq 3A> 3$ for $r\geq S$.

Choosing 
$$R:=r\left(1+\frac{1}{(\log T(r))^2}\right)$$
in (\ref{MT}) and noting that $\log(T/S)> 2>\pi^2/6$
we deduce from Lemma~\ref{borel} that there exists $r\in [S,T]$
such that 
$$
\log M(r,h)\leq 
\left(2+\frac{1}{(\log T(r,h))^2}\right) (\log T(r,h))^2 eT(r,h).$$
Now $(\log T(r,h))^2\geq (\log 3A)^2\geq 1$ and hence
\begin{equation}\label{MT2}
\log M(r,h)
\leq 
3eT(r,h) (\log T(r,h))^2.
\end{equation}
For a value of $r$ satisfying (\ref{MT2}) we consider the set
$$E(r):=\left\{t\in [0,2\pi ]: \log|h(re^{it})|\geq\frac12 T(r,h)\right\}$$
and define
$\lambda(r):=\meas E(r)$, where $\meas E$ denotes the measure
of a set $E$.
Then for at least one the four sets $E_\ell(r):=E(r)\cap [\ell\frac{\pi}{2},
(\ell+1)\frac{\pi}{2} ]$,  where $\ell\in\{0,1,2,3\}$,  we have 
$\meas E_\ell(r)\geq \frac14\lambda(r)$.
We assume now that this is the case for $\ell=0$. The modifications 
that have to be made for the other cases will be obvious.
We define 
$$\alpha(r):=\max E_0(r)$$
and
$$\beta(r):=\min\left\{t\in  \left[\alpha(r),\frac{\pi}{2}\right]:
\log|h(re^{it})|=A\right\}$$
Then $\alpha(r)\geq \meas E_0(r)\geq \frac14\lambda(r)$.
For $\alpha(r)\leq t\leq \beta(r)$ we have 
$1<A\leq  \log|h(re^{it})|\leq\frac12 T(r,h)$.
Thus 
$\log\log h$ may be defined on the arc
$\{re^{it}:\alpha(r)\leq t\leq \beta(r)\}$, and we may choose the 
branch of the logarithm such that 
$|\log h(re^{i\beta(r)})|\leq \log |h(re^{i\beta(r)})|+\pi=A+\pi$
and hence 
$\log |\log h(re^{i\beta(r)})|\leq \log (A+\pi)$.
Thus
\begin{eqnarray*}
\log \log| h(re^{i\alpha(r)})|- \log (A+\pi)
&\leq &
\log |\log h(re^{i\alpha(r)})|- \log |\log h(re^{i\beta(r)})|\\ 
&\leq &
|\log \log h(re^{i\alpha(r)})- \log \log h(re^{i\beta(r)})|\\ 
&=&
\left|
\int_{\alpha(r)}^{\beta(r)} 
\frac{h'(re^{it})|}{h(re^{it})\log h(re^{it})}ire^{it}dt\right|\\
&\leq &
r \int_{\alpha(r)}^{\beta(r)} 
\frac{|h'(re^{it})|}{|h(re^{it})|\log|h(re^{it})|}dt\\
&\leq &
4 \int_{\alpha(r)}^{\beta(r)} 
\frac{dt}{|\sin t|}.
\end{eqnarray*}
Here the last inequality follows from (\ref{hlogh}).
Now $\sin t\geq 2t/\pi$ for $0\leq t\leq \pi/2$. Thus
$$\log \log| h(re^{i\alpha(r)})|- \log (A+\pi)
\leq 2\pi \int_{\alpha(r)}^{\beta(r)} 
\frac{dt}{t}=
 2\pi \log \frac{\beta(r)}{\alpha(r)} 
 \leq 
 2\pi \log \frac{\pi}{2\alpha(r)} .$$
Since 
$$\log|h(re^{i\alpha(r)})|=\frac12 T(r,h)$$
this yields 
\begin{equation}\label{logT}
\log T(r,h)\leq 
 2\pi \log \frac{\pi}{2\alpha(r)} 
 +\log(A+\pi)+\log 2.
\end{equation}
 Now 
\begin{eqnarray*}
T(r,h)
&=&
\frac{1}{2\pi}\int_0^{2\pi} \log^+ |h(re^{it})| dt\\
&= &
\frac{1}{2\pi}\int_{E(r)} \log^+ |h(re^{it})| dt
+
\frac{1}{2\pi}\int_{[0,2\pi]\backslash E(r)} \log^+ |h(re^{it})| dt\\
&\leq &
\frac{1}{2\pi}\lambda(r)\log M(r,h) +
\frac12 T(r,h)
\end{eqnarray*}
so that 
$$\frac{\lambda(r)}{\pi}\geq \frac{T(r,h)}{\log M(r,h)}.$$
Hence
$$ \frac{\pi}{2\alpha(r)}\leq 
\frac{2\pi}{\lambda(r)}\leq 
\frac{2\log M(r,h)}{T(r,h)}.$$
Using (\ref{MT2}) we find that
$$ \frac{\pi}{2\alpha(r)}\leq 
6 e (\log T(r,h))^2.$$
Together with (\ref{logT}) this yields 
$$\log T(r,h)\leq 
 2\pi \log 
\left(6 e (\log T(r,h))^2\right) 
 +\log(A+\pi)+\log 2.$$
 This implies that 
$$\log T(r,h)\leq 
4\pi\log\log T(r,h)+C$$
for some constant $C$, which is a contradiction
since $T(r,h)=T(r,g_n)\to\infty$ as $n\to\infty$.\qed

\noindent 

\end{document}